# Bipolar Neutrosophic Soft Expert Sets

Mehmet Şahin, İrfan Deli and Vakkas Uluçay

Department of Mathematics, Gaziantep University, Gaziantep27310-Turkey
Email: mesahin@gantep.edu.tr, vulucay27@gmail.com

Muallim Rıfat Faculty of Education, Kilis 7 Aralık University, 79000 Kilis, Turkey,
irfandeli@kilis.edu.tr

**Abstract**

In this paper, we introduce concept of bipolar neutrosophic soft expert set and its some operations. Also, we propose score, certainty and accuracy functions to compare the bipolar neutrosophic soft expert sets. We give examples for these concepts.

## 1. Introduction

In some real life problems in expert system, belief system, information fusion and so on, we must consider the truth-membership as well as the falsity- membership for proper description of an object in uncertain, ambiguous environment. Intuitionistic fuzzy sets introduced by Atanassov [1]. After Atanassov's work, Smarandache [17] introduced the concept of neutrosophic set which is a mathematical tool for handling problems involving imprecise, indeterminacy and inconsistent data. These sets models have been studied by many authors; on application [4-7,10-12,15,16], and so on.

Bosc and Pivert [2] said that "Bipolarity refers to the propensity of the human mind to reason and make decisions on the basis of positive and negative effects. Positive information states what is possible, satisfactory, permitted, desired, or considered as being acceptable. On the other hand, negative statements express what is impossible, rejected, or forbidden. Negative preferences correspond to constraints, since they specify which values or objects have to be rejected (i.e., those that do not satisfy the constraints), while positive preferences correspond to wishes, as they specify which objects are more desirable than others (i.e., satisfy user wishes) without rejecting those that do not meet the wishes." Therefore, Lee [8,9] introduced the concept of bipolar fuzzy sets which is an generalization of the fuzzy sets. Recently, bipolar fuzzy models have been studied by many authors on algebraic structures such as; Majumder [13] proposed bipolar valued fuzzy subsemigroup, bipolar valued fuzzy bi-ideal, bipolar valued fuzzy (1, 2) - ideal and bipolar valued fuzzy ideal. Manemaran and Chellappa [14] gave some applications of bipolar fuzzy sets in groups are called the bipolar fuzzy groups, fuzzy d-ideals of groups under (T-S) norm. Chen et al. [3] studied of $m$-polar fuzzy set and illustrates how many concepts have been defined based on bipolar fuzzy sets.

Alkhazaleh et al. [21] where the mapping in which the approximate function are defined from fuzzy parameters set, and gave an application of this concept in decision making. Alkhazaleh and Salleh [22] introduced the concept soft expert sets where user can know the opinion of all expert sets. Sahin et al. [23] firstly proposed neutrosophic soft expert sets with operations.

In this paper, we introduced the concept of bipolar neutrosophic soft expert sets which is an extension of the fuzzy soft expert sets, bipolar fuzzy sets, intuitionistic fuzzy sets and neutrosophic sets. Also, we give some operations and operators on the bipolar neutrosophic soft expert sets. In section1, from intuitionistic fuzzy sets to bipolar neutrosophic sets are mentions. In section2, preliminaries are given. In section3, the concept of bipolar neutrosophic soft expert set and its basic operations, namely complement, union and intersection. In Section 4 give conclusions.

## 2. Preliminaries

In this section we recall some related definitions.

**2.1. Definition:** [17] Let U be a space of points (objects), with a generic element in U denoted by u. A neutrosophic sets (N-sets) A in U is characterized by a truth-membership function $T_A$, a indeterminacy-membership function $I_A$ and a falsity-membership function $F_A$. $T_A$ (u); $I_A$ (u) and $F_A$ (u) are real standard or nonstandard subsets of [0, 1]. It can be written as
A = {< u, ($T_A$ (u), $I_A$ (u), $F_A$ (u)) >: u ∈ U, $T_A$ (u), $I_A$ (u), $F_A$ (u) ∈ [0, 1]}. There is no restriction on the sum of $T_A$ (u); $I_A$ (u) and $F_A$ (u), so
$$0 \leq \sup T_A (u) + \sup I_A (u) + \sup F_A (u) \leq 3.$$

**2.2. Definition:** [20] A neutrosophic set A is contained in another neutrosophic set B i.e. $A \subseteq B$ if $\forall x \in X$, $T_A(x) \leq T_B(x)$, $I_A(x) \leq I_B(x)$, $F_A(x) \geq F_B(x)$.

Let U be a universe, E a set of parameters, and X a soft experts (agents). Let $O$ be a set of opinion, $Z = E \times X \times O$ and $A \subseteq Z$.

**2.3. Definition:** [23] A pair $(F, A)$ is called a neutrosophic soft expert set over U, where F is mapping given by

$$F: A \to P(U)$$

Where $P(U)$ denotes the power neutrosophic set of U.

Set- theoretic operations, for two neutrosophic soft expert sets,

$A_{NSE}$ = {<x, $T_{F(e)}$(x), $I_{F(e)}$(x), $F_{F(e)}$(x)> | $\forall e \in A, x \in U$ } and $B_{NSE}$ = {<x, $T_{G(e)}$(x), $I_{G(e)}$(x), $F_{G(e)}$(x)> | $\forall e \in A, x \in U$ } are given as;

1. The subset; $A_{NSE} \subseteq B_{NSE}$ if and only if

$$T_{F(e)}(x) \widetilde{\leq} T_{G(e)}(x), \; I_{F(e)}(x) \widetilde{\leq} I_{G(e)}(x), \; F_{F(e)}(x) \widetilde{\geq} F_{G(e)}(x) \; \forall e \in A, x \in U.$$

2. $A_{NSE} = B_{NSE}$ if and only if,

$$T_{F(e)}(x) = T_{G(e)}(x), \; I_{F(e)}(x) = I_{G(e)}(x), \; F_{F(e)}(x) = F_{G(e)}(x) \; \forall e \in A, x \in U.$$

3. The complement of $A_{NSE}$ is denoted by $A_{NS}^c$ and is defined by

$$A_{NSE}^c = \{<x, \; T_{F^c(x)} = F_{F(x)}, \; I_{F^c(x)} = I_{F(x)}, \; F_{F^c(x)} = T_{F(x)} | \; x \in X \}$$

4. The intersection

$$A_{NSE} \cap B_{NSE} = \{<x, \min\{T_{F(e)}(x), T_{G(e)}(x)\}, \max\{I_{F(e)}(x), I_{G(e)}(x)\},$$
$$\max\{F_{F(e)}(x), F_{G(e)}(x)\}>: x \in X \}$$

5. The union

$$A_{NSE} \cup B_{NSE} = \{<x, \max\{T_{F(e)}(x), T_{G(e)}(x)\}, \min\{I_{F(e)}(x), I_{G(e)}(x)\},$$
$$\min\{F_{F(e)}(x), F_{G(e)}(x)\}>: x \in X \}$$

**2.4. Definition:** [27] A bipolar neutrosophic set $A$ in $X$ is defined as an object of the form
$$A = \{\langle x, T^+(x), I^+(x), F^+(x), T^-(x), I^-(x), F^-(x)\rangle : x \in X\},$$
where $T^+, I^+, F^+ : X \to [1,0]$ and $T^-, I^-, F^- : X \to [-1,0]$.

**2.5. Definition:** [27] Let $\tilde{a}_1 = \langle T_1^+, I_1^+, F_1^+, T_1^-, I_1^-, F_1^- \rangle$ and $\tilde{a}_2 = \langle T_2^+, I_2^+, F_2^+, T_2^-, I_2^-, F_2^- \rangle$ be two bipolar neutrosophic number. Then the operations for NNs are defined as below;

i. $\lambda \tilde{a}_1 = \langle 1 - (1 - T_1^+)^\lambda, (I_1^+)^\lambda, (F_1^+)^\lambda, -(-T_1^-)^\lambda, -(-I_1^-)^\lambda, -(1 - (1 - (-F_1^-))^\lambda) \rangle$

ii. $\tilde{a}_1^\lambda =$
$\langle (T_1^+)^\lambda, 1 - (1 - I_1^+)^\lambda, 1 - (1 - F_1^+)^\lambda, -(1 - (1 - (-T_1^-))^\lambda), -(-I_1^-)^\lambda, -(-F_1^-)^\lambda \rangle$

iii. $\tilde{a}_1 + \tilde{a}_2 =$
$\langle T_1^+ + T_2^+ - T_1^+ T_2^+, I_1^+ I_2^+, F_1^+ F_2^+, -T_1^- T_2^-, -(-I_1^- - I_2^- - I_1^- I_2^-), -(-F_1^- - F_2^- - F_1^- F_2^-) \rangle$

iv. $\tilde{a}_1 \cdot \tilde{a}_2 =$
$\langle T_1^+ T_2^+, I_1^+ + I_2^+ - I_1^+ I_2^+, F_1^+ + F_2^+ - F_1^+ F_2^+, -(-T_1^- - T_2^- - T_1^- T_2^-), -I_1^- I_2^-, -F_1^- F_2^- \rangle$

where $\lambda > 0$.

**2.6. Definition:** [27] Let $\tilde{a}_1 = \langle T_1^+, I_1^+, F_1^+, T_1^-, I_1^-, F_1^- \rangle$ be a bipolar neutrosophic number. Then, the score function $s(\tilde{a}_1)$, accuracy function $a(\tilde{a}_1)$ and certainty function $c(\tilde{a}_1)$ of an NBN are defined as follows:

i. $\tilde{s}(\tilde{a}_1) = (T_1^+ + 1 - I_1^+ + 1 - F_1^+ + 1 + T_1^- - I_1^- - F_1^-)/6$
ii. $\tilde{a}(\tilde{a}_1) = T_1^+ - F_1^+ + T_1^- - F_1^-$
iii. $\tilde{c}(\tilde{a}_1) = T_1^+ - F_1^-$

**2.7. Definition:** [27] $\tilde{a}_1 = \langle T_1^+, I_1^+, F_1^+, T_1^-, I_1^-, F_1^- \rangle$ and $\tilde{a}_2 = \langle T_2^+, I_2^+, F_2^+, T_2^-, I_2^-, F_2^- \rangle$ be two bipolar neutrosophic number. The comparison method can be defined as follows:

i. if $\tilde{s}(\tilde{a}_1) > \tilde{s}(\tilde{a}_2)$, then $\tilde{a}_1$ is greater than $\tilde{a}_2$, that is, $\tilde{a}_1$ is superior to $\tilde{a}_2$, denoted by $a_1 > \tilde{a}_2$

ii. $\tilde{s}(\tilde{a}_1) = \tilde{s}(\tilde{a}_2)$ and $\tilde{a}(\tilde{a}_1) > \tilde{a}(\tilde{a}_2)$, then $\tilde{a}_1$ is greater than $\tilde{a}_2$, that is, $\tilde{a}_1$ is superior to $\tilde{a}_2$, denoted by $\tilde{a}_1 < \tilde{a}_2$;

iii. if $\tilde{s}(\tilde{a}_1) = \tilde{s}(\tilde{a}_2)$, $\tilde{a}(\tilde{a}_1) = \tilde{a}(\tilde{a}_1)$ and $\tilde{c}(\tilde{a}_1) > \tilde{c}(\tilde{a}_2)$, then $\tilde{a}_1$ is greater than $\tilde{a}_2$, that is, $\tilde{a}_1$ is superior to $\tilde{a}_2$, denoted by $\tilde{a}_1 > \tilde{a}_2$;

iv. if $\tilde{s}(\tilde{a}_1) = \tilde{s}(\tilde{a}_2)$, $\tilde{a}(\tilde{a}_1) = \tilde{a}(\tilde{a}_2))$ and $\tilde{c}(\tilde{a}_1) = \tilde{c}(\tilde{a}_2)$, then $\tilde{a}_1$ is equal to $\tilde{a}_2$, that is, $\tilde{a}_1$ is indifferent to $\tilde{a}_2$, denoted by $\tilde{a}_1 = \tilde{a}_2$.

## 3. Bipolar Neutrosophic Soft Expert Set

In this section, using the concept of bipolar neutrosophic set now we introduce the concept of bipolar neutrosophic soft expert set and we also give basic properties of this concept.

Let U be a universe, E a set of parameters, X a set of experts (agents), and $O = \{1 = agree, 0 = disagree\}$ a set of opinions. Let $Z = E \times X \times O$ and $\bar{A} \subseteq Z$.

**3.1. Definition:** A pair $(H, \bar{A})$ is called a bipolar neutrosophic soft expert set over U, where H is mapping given by

$$H: \bar{A} \to P(U)$$

where $P(U)$ denotes the power bipolar neutrosophic set of U and

$$(H, \bar{A}) = \{\langle u, T^+_{H(e)}(u), I^+_{H(e)}(u), F^+_{H(e)}(u), T^-_{H(e)}(u), I^-_{H(e)}(u), F^-_{H(e)}(u) \rangle : \forall e \in A, u \in U\},$$

where $T^+_{H(e)}, I^+_{H(e)}, F^+_{H(e)} : U \to [1, 0]$ and $T^-_{H(e)}, I^-_{H(e)}, F^-_{H(e)} : U \to [-1, 0]$.

For definition we consider an example.

**3.2. Example:** Suppose the following $U$ is the set of notebook under consideration $E$ is the set of parameters. Each parameter is a neutrosophic word or sentence involving neutrosophic words.

$E = \{\text{cheap}; \text{expensive}\} = \{e_1, e_2\}$

$X = \{p, q, r\}$ be a set of experts. Suppose that

$H(e_1, p, 1) = \{< u_1, 0.3, 0.5, 0.7, -0.2, -0.3, -0.4 >, < u_3, 0.5, 0.6, 0.3, -0.3, -0.4, -0.1 >\}$

$H(e_1, q, 1) = \{< u_2, 0.8, 0.2, 0.3, -0.1, -0.3, -0.5 >, < u_3, 0.9, 0.5, 0.7, -0.4, -0.1, -0.2 >\}$

$H(e_1, r, 1) = \{< u_1, 0.4, 0.7, 0.6, -0.6, -0.2, -0.4 >\}$

$H(e_2, p, 1) = \{< u_1, 0.4, 0.2, 0.3, -0.2, -0.3, -0.1 >, < u_2, 0.7, 0.1, 0.3, -0.3, -0.2, -0.5 >\}$

$H(e_2, q, 1) = \{< u_3, 0.3, 0.4, 0.2, -0.5, -0.1, -0.4 >\}$

$H(e_2, r, 1) = \{< u_2, 0.3, 0.4, 0.9, -0.4, -0.3, -0.1 >\}$

$H(e_1, p, 0) = \{< u_2, 0.5, 0.2, 0.3, -0.5, -0.2, -0.3 >\}$

$H(e_1, q, 0) = \{< u_1, 0.6, 0.3, 0.5, -0.4, -0.2, -0.6 >\}$

$H(e_1, r, 0) = \{< u_2, 0.7, 0.6, 0.4, -0.3, -0.4, -0.5 >, < u_3, 0.9, 0.5, 0.7, -0.2, -0.3, -0.5 >\}$

$H(e_2, p, 0) = \{< u_3, 0.7, 0.9, 0.6, -0.2, -0.3, -0.4 >\}$

$H(e_2, q, 0) = \{< u_1, 0.7, 0.3, 0.6, -0.3, -0.2, -0.4 >, < u_2, 0.6, 0.2, 0.5, -0.3, -0.1, -0.4 >\}$

$H(e_2, r, 0) = \{< u_1, 0.6, 0.2, 0.5, -0.5, -0.3, -0.2 >, < u_3, 0.7, 0.2, 0.8, -0.6, -0.2, -0.1 >\}$

The bipolar neutrosophic soft expert set $(H, \bar{Z})$ is a parameterized family $\{H(e_i), i = 1, 2, 3, ...\}$ of all neutrosophic sets of $U$ and describes a collection of approximation of an object.

**3.3. Definition:** Let $(H, \bar{A})$ and $(G, \bar{B})$ be two bipolar neutrosophic soft expert sets over the common universe U. $(H, \bar{A})$ is said to be bipolar neutrosophic soft expert subset of $(G, \bar{B})$, if $(H, \bar{A}) \subseteq (G, \bar{B})$ if and only if

$$T^+_{H(e)}(u) \leq T^+_{G(e)}(u) \; I^+_{H(e)}(u) \leq I^+_{G(e)}(u), F^+_{H(e)}(u) \geq F^+_{G(e)}(u),$$

and

$$T^-_{H(e)}(u) \geq T^-_{G(e)}(u), \; I^-_{H(e)}(u) \geq I^-_{G(e)}(u), \; F^-_{H(e)}(u) \leq F^-_{G(e)}(u)$$

$\forall e \in A, u \in U.$

$(H, \bar{A})$ is said to be bipolar neutrosophic soft expert superset of $(G, \bar{B})$ if $(G, \bar{B})$ is a neutrosophic soft expert subset of $(H, \bar{A})$. We denote by $(H, \bar{A}) \tilde{\supseteq} (G, \bar{B})$.

**3.4 Example**: Suppose that a company produced new types of its products and wishes to take the opinion of some experts about price of these products. Let $U = \{u_1, u_2, u_3\}$ be a set of product, $E = \{e_1, e_2\}$ a set of decision parameters where $e_i (i = 1,2)$ denotes the decision "cheap ", "expensive" respectively and let $X = \{p, q, r\}$ be a set of experts. Suppose $(H, \bar{A})$ and $(G, \bar{B})$ be defined as follows:

$(H, \bar{A}) =$

$\{[(e_1, p, 1), < u_1, 0.3, 0.5, 0.6, -0.2, -0.3, -0.4 >, < u_2, 0.5, 0.2, 0.3, -0.4, -0.2, -0.5 >],$

$[(e_2, p, 0), < u_2, 0.2, 0.4, 0.7, -0.5, -0.4, -0.3 >],$
$[(e_1, q, 1), < u_1, 0.6, 0.3, 0.5, -0.6, -0.2, -0.5 >, < u_2, 0.6, 0.2, 0.3, -0.5, -0.4, -0.3 >],$

$[(e_1, r, 0), < u_1, 0.2, 0.7, 0.3, -0.4, -0.3, -0.5 >],$
$[(e_2, r, 1), < u_2, 0.3, 0.4, 0.9, -0.3, -0.2, -0.4 >, < u_3, 0.7, 0.2, 0.8, -0.5, -0.3, -0.6 >]\}.$

$(G, \bar{B}) =$

$\{[(e_1, p, 1), < u_1, 0.3, 0.5, 0.7, -0.2, -0.3, -0.6 >, < u_2, 0.5, 0.2, 0.3, -0.1, -0.2, -0.7 >],$

$[(e_2, p, 0), < u_2, 0.2, 0.4, 0.7, -0.2, -0.4, -0.5 >],$

$[(e_1, q, 1), < u_1, 0.6, 0.3, 0.5, -0.1, -0.2, -0.8 >, < u_2, 0.6, 0.2, 0.3, -0.3 - 0.1, -0.4 >]\}.$

Therefore

$$(H, \bar{A}) \tilde{\supseteq} (G, \bar{B}).$$

**3.5. Definition:** Let $(H, \bar{A})$ and $(G, \bar{B})$ be two bipolar neutrosophic soft expert sets over the common universe U. $(H, \bar{A})$ is said to be bipolar neutrosophic soft expert equal $(G, B)$, if $(H, \bar{A}) = (G, \bar{B})$ if and only if

$$T^+_{H(e)}(u) = T^+_{G(e)}(u) \; I^+_{H(e)}(u) = I^+_{G(e)}(u), F^+_{H(e)}(u) = F^+_{G(e)}(u),$$

and

$$T^-_{H(e)}(u) = T^-_{G(e)}(u), \; I^-_{H(e)}(u) = I^-_{G(e)}(u), \; F^-_{H(e)}(u) = F^-_{G(e)}(u)$$

$\forall e \in A, u \in U.$

**3.6. Definition:** NOT set of set parameters. Let $E = \{e_1, e_2, \ldots, e_n\}$ be a set of parameters. The NOT set of E is denoted by $\neg E = \{\neg e_1, \neg e_2, \ldots, \neg e_n\}$ where $\neg e_i$ = not $e_i$, $\forall$ i=1,2,...,n.

**3.7. Example:** Consider 3.2 example. Here ¬E={not cheap, not expensive}

**3.8. Definition:** Complement of a bipolar neutrosophic soft expert set. The complement of a bipolar neutrosophic soft expert set $(H, \bar{A})$ denoted by $(H, \bar{A})^c$ and is defined as $(H, \bar{A})^c = (H^c, \neg \bar{A})$ where $H^c = \neg \bar{A} \to P(U)$ is mapping given by $H^c(u)$= neutrosophic soft expert complement with

$$T^+_{H^c(u)} = F^+_{H(u)}, I^+_{H^c(u)} = I^+_{H(u)}, F^+_{H^c(u)} = T^+_{H(u)}$$

and

$$T^-_{H^c(u)} = F^-_{H(u)}, I^-_{H^c(u)} = I^-_{H(u)}, F^-_{H^c(u)} = T^-_{H(u)}$$

**3.9. Example:** Consider the 3.2 Example. Then $(H, \bar{Z})^c$ describes the "not price of the notebook" we have

$(H, \bar{Z})^c = \{(\neg e_1, p, 1), [< u_2, 0.3, 0.2, 0.5, -0.3, -0.2, -0.5 >]$

$[(\neg e_1, q, 1), < u_1, 0.5, 0.3, 0.6, -0.4, -0.1, -0.3 >],$

$[(\neg e_1, r, 1), < u_2, 0.4, 0.6, 0.7, -0.3, -0.4, -0.2 >, < u_3, 0.7, 0.5, 0.9, -0.1, -0.2, -0.3 >],$

$[(\neg e_2, p, 1), < u_3, 0.6, 0.9, 0.7, -0.4, -0.3, -0.2 >],$

$[(\neg e_2, q, 1), < u_1, 0.6, 0.3, 0.7, -0.5, -0.1, -0.3 >, < u_2, 0.5, 0.2, 0.6, -0.3, -0.5, -0.6 >],$

$[(\neg e_2, r, 1), < u_1, 0.5, 0.2, 0.6 - 0.6, -0.2, -0.4 >, < u_3, 0.8, 0.2, 0.7, -0.3, -0.4, -0.1 >],$

$[(\neg e_1, p, 0), < u_1, 0.7, 0.5, 0.3, -0.4, -0.2, -0.3 >, < u_3, 0.3, 0.6, 0.5, -0.6, -0.3, -0.5 >],$

$[(\neg e_1, q, 0), < u_2, 03, 0.2, 0.8, -0.3, -0.2, -0.7 >, < u_3, 0.9, 0.5, 0.7, -0.7, -0.3, -0.5 >],$

$[(\neg e_1, r, 0), < u_1, 0.6, 0.7, 0.4, -0.4, -0.3, -0.5 >],$

$[(\neg e_2, p, 0), < u_1, 0.3, 0.2, 0.4, -0.3, -0.5, -0.4 >, < u_2, 0.3, 0.1, 0.7, -0.6, -0.5, -0.1 >],$

$[(\neg e_2, q, 0), < u_3, 0.2, 0.4, 0.3, -0.7, -0.4, -0.3 >],$

$[(\neg e_2, r, 0), < u_2, 0.9, 0.4, 0.3, -0.8, -0.3, -0.5 >]\}.$

**3.10 Definition:** Empty or Null bipolar neutrosophic soft expert set with respect to parameter. A bipolar neutrosophic soft expert set $(H, \bar{A})$ over the universe $U$ is termed to be empty or null bipolar neutrosophic soft expert set with respect to the parameter $\bar{A}$ if

$$T^+_{H(e)}(u) = T^+_{G(e)}(u) = 0 \; I^+_{H(e)}(u) = I^+_{G(e)}(u) = 0, F^+_{H(e)}(u) = F^+_{G(e)}(u) = 0,$$

and

$$T^-_{H(e)}(u) = T^-_{G(e)}(u) = 0, \; I^-_{H(e)}(u) = I^-_{G(e)}(u) = 0, \; F^-_{H(e)}(u) = F^-_{G(e)}(u) = 0$$

$\forall e \in \bar{A}, u \in U$.

In this case the null bipolar neutrosophic soft expert set (NBNSES) is denoted by $\phi_{\breve{A}}$.

**3.11 Example:** Let $U = \{u_1, u_2, u_3\}$ the set of three handbags be considered as universal set $E = \{quality\} = \{e_1\}$ be the set of parameters that characterizes the handbag and let $X = \{p, q\}$ be a set of experts.

$\Phi_{\breve{A}} = (NBNSES) = \{[(e_1, p, 1), < u_1, 0,0,0,0,0,0 >, < u_2, 0,0,0,0,0,0 >],$

$$[(e_1, q, 1), < u_1, 0,0,0,0,0,0 >, < u_2, 0,0,0,0,0,0 >],$$
$$[(e_1, p, 0), < u_3, 0,0,0,0,0,0 >],$$
$$[(e_1, q, 0), < u_3, 0,0,0,0,0,0 >]\}.$$

Here the (NBNSES) (H, $\bar{A}$) is the null bipolar neutrosophic soft expert sets.

**3.12 Definition:** An agree-bipolar neutrosophic soft expert set $(H, \bar{A})_1$ over $U$ is a bipolar neutrosophic soft expert subset of $(H, \bar{A})$ defined as follow

$$(H, \bar{A})_1 = \{H_1(u): u \in E \times X \times \{1\}\}.$$

**3.13 Example:** Consider 3.2. Example. Then the agree-bipolar neutrosophic soft expert set $(H, \bar{A})_1$ over $U$ is
$(H, \bar{A})_1 =$
$\{[(e_1, p, 1), < u_1, 0.3, 0.5, 0.7, -0.2, -0.3, -0.4 >, < u_3, 0.5, 0.6, 0.3, -0.3, -0.4, -0.1 >],$

$[(e_1, q, 1), < u_2, 0.8, 0.2, 0.3, -0.1, -0.3, -0.5 >, < u_3, 0.9, 0.5, 0.7, -0.4, -0.1, -0.2 >],$

$[(e_1, r, 1), < u_1, 0.4, 0.7, 0.6, -0.6, -0.2, -0.4 >],$

$[(e_2, p, 1), < u_1, 0.4, 0.2, 0.3, -0.2, -0.3, -0.1 >, < u_2, 0.7, 0.1, 0.3, -0.3, -0.2, -0.5 >],$

$[(e_2, q, 1), < u_3, 0.3, 0.4, 0.2, -0.5, -0.1, -0.4 >],$

$[(e_2, r, 1), < u_2, 0.3, 0.4, 0.9, -0.4, -0.3, -0.1 >]\}.$

**3.14 Definition:** A disagree-bipolar neutrosophic soft expert set $(H, \bar{A})_0$ over U is a bipolar neutrosophic soft expert subset of $(H, \bar{A})$ defined as follow

$$(H, \bar{A})_0 = \{F_0(u): u \in E \times X \times \{0\}\}.$$

**3.15 Example:** Consider 3.2 Example. Then the disagree-bipolar neutrosophic soft expert set $(H, \bar{A})_0$ over $U$ is
$(H, \bar{A})_0 = \{[(e_1, p, 0), < u_2, 0.5, 0.2, 0.3, -0.5, -0.2, -0.3 >],$

$[(e_1, q, 0), < u_1, 0.6, 0.3, 0.5, -0.4, -0.2, -0.6 >],$

$[(e_1, r, 0), < u_2, 0.7, 0.6, 0.4, -0.3, -0.4, -0.5 >, < u_3, 0.9, 0.5, 0.7, -0.2, -0.3, -0.5 >],$

$[(e_2, p, 0), < u_3, 0.7, 0.9, 0.6, -0.2, -0.3, -0.4 >],$

$[(e_2, q, 0), < u_1, 0.7, 0.3, 0.6, -0.3, -0.2, -0.4 >, < u_2, 0.6, 0.2, 0.5, -0.3, -0.1, -0.4 >],$

$[(e_2, r, 0), < u_1, 0.6, 0.2, 0.5, -0.5, -0.3, -0.2 >, < u_3, 0.7, 0.2, 0.8, -0.6, -0.2, -0.1 >]\}.$

**3.16 Definition:** Union of two bipolar neutrosophic soft expert sets. Let

$(H, \bar{A}) = \{\langle u, T^+_{H(e)}(u), I^+_{H(e)}(u), F^+_{H(e)}(u), T^-_{H(e)}(u), I^-_{H(e)}(u), F^-_{H(e)}(u) \rangle : \forall e \in A, u \in U\}$ and

$(G, \bar{B}) = \{\langle u, T^+_{G(e)}(u), I^+_{G(e)}(u), F^+_{G(e)}(u), T^-_{G(e)}(u), I^-_{G(e)}(u), F^-_{G(e)}(u) \rangle : \forall e \in B, u \in U\}$ be two bipolar neutrosophic soft expert sets. Then their union is defined as:

$$((H,\bar{A})\tilde{\cup}(G,\bar{B}))(u) = \begin{pmatrix} \max(T^+_{H(e)}(u), T^+_{G(e)}(u)), \dfrac{I^+_{H(e)}(u)+I^+_{G(e)}(u)}{2}, \min((F^+_{H(e)}(u), F^+_{G(e)}(u)), \\ \min(T^-_{H(e)}(u), T^-_{G(e)}(u)), \dfrac{I^-_{H(e)}(u)+I^-_{G(e)}(u)}{2}, \max((F^-_{H(e)}(u), F^-_{G(e)}(u)) \end{pmatrix}$$

$\forall e \in A, u \in U.$

**3.17 Example:** Let $(H, \bar{A})$ and $(G, \bar{B})$ be two BNSESs over the common universe $U$

$(H, \bar{A}) =$
$\{[(e_1, p, 1), < u_1, 0.2, 0.5, 0.8, -0.4, -0.3, -0.5 >, < u_3, 0.2, 0.6, 0.5, -0.2, -0.1, -0.4 >],$
$[(e_1, q, 1), < u_1, 0.5, 0.3, 0.6, -0.2, -0.1, -0.3 >, < u_2, 0.8, 0.2, 0.3, -0.2, -0.3, -0.1 >]\}$

$(G, \bar{B})$
$= \{(e_1, p, 1), < u_1, 0.1, 0.6, 0.2, -0.3, -0.1, -0.4 >, < u_2, 0.4, 0.5, 0.8, -0.1, -0.3, -0.5 >\}$

Therefore $(H, \bar{A}) \tilde{\cup} (G, \bar{B}) = (R, \bar{C})$

$(R, \bar{C}) = \left\{ \begin{bmatrix} (e_1, p, 1), < u_1, 0.2, 0.55, 0.2, -0.4, -0.2, -0.4 >, \\ < u_2, 0.4, 0.5, 0.8 - 0.1, -0.3, -0.5 >, < u_3, 0.2, 0.6, 0.5, -0.2, -0.1, -0.4 > \end{bmatrix}, \right.$

$[(e_1, q, 1), < u_1, 0.5, 0.3, 0.6, -0.2, -0.1, -0.3 >, < u_2, 0.4, 0.5, 0.8, -0.1, -0.3, -0,5 >]\}.$

**3.18 Definition:** Intersection of two bipolar neutrosophic soft expert sets.

$(H, \bar{A}) = \{\langle u, T^+_{H(e)}(u), I^+_{H(e)}(u), F^+_{H(e)}(u), T^-_{H(e)}(u), I^-_{H(e)}(u), F^-_{H(e)}(u) \rangle : \forall e \in A, u \in U\}$ and

$(G, \bar{B}) = \{\langle u, T^+_{G(e)}(u), I^+_{G(e)}(u), F^+_{G(e)}(u), T^-_{G(e)}(u), I^-_{G(e)}(u), F^-_{G(e)}(u) \rangle : \forall e \in B, u \in U\}$ be two

bipolar neutrosophic soft expert sets. Then their intersection is defined as:

$$((H,\bar{A})\tilde{\cap}(G,\bar{B}))(u) = \begin{pmatrix} \min(T^+_{H(e)}(u), T^+_{G(e)}(u)), \dfrac{I^+_{H(e)}(u)+I^+_{G(e)}(u)}{2}, \max((F^+_{H(e)}(u), F^+_{G(e)}(u)), \\ \max(T^-_{H(e)}(u), T^-_{G(e)}(u)), \dfrac{I^-_{H(e)}(u)+I^-_{G(e)}(u)}{2}, \min((F^-_{H(e)}(u), F^-_{G(e)}(u)) \end{pmatrix}$$

$\forall e \in A, u \in U.$

**3.19 Example: :** Let $(H, \bar{A})$ and $(G, \bar{B})$ be two BNSESs over the common universe $U$

$(H, \bar{A}) =$
$\{[(e_1, p, 1), < u_1, 0.2, 0.5, 0.8, -0.4, -0.3, -0.5 >, < u_3, 0.2, 0.6, 0.5, -0.2, -0.1, -0.4 >],$
$[(e_1, q, 1), < u_1, 0.5, 0.3, 0.6, -0.2, -0.1, -0.3 >, < u_2, 0.8, 0.2, 0.3, -0.2, -0.3, -0.1 >]\}$

$(G, \bar{B})$
$= \{(e_1, p, 1), < u_1, 0.1, 0.6, 0.2, -0.3, -0.1, -0.4 >, < u_2, 0.4, 0.5, 0.8, -0.1, -0.3, -0.5 >\}$

Therefore $(H, \bar{A}) \tilde{\cap} (G, \bar{B}) = (R, \bar{C})$

$(R, \bar{C}) = \{[(e_1, p, 1), < u_1, 0.1, 0.55, 0.8, -0.3, -0.2, -0.5 >]\}.$

**3.1. Proposition:** If $(H, \bar{A})$ and $(G, B)$ are bipolar neutrosophic soft expert sets over $U$. Then

i. $(H,\bar{A}) \,\tilde{\cup}\, (G,\bar{B}) = (G,\bar{B}) \,\tilde{\cup}\, (H,\bar{A})$

ii. $(H,\bar{A}) \,\tilde{\cap}\, (G,\bar{B}) = (G,\bar{B}) \,\tilde{\cap}\, (H,\bar{A})$

iii. $((H,\bar{A})^c)^c = (H,\bar{A})$

iv. $(H,\bar{A}) \,\tilde{\cup}\, \phi = (H,\bar{A}),\qquad (H,\bar{A}) \,\tilde{\cap}\, \phi = \phi$

**Proof:** The proof is straightforward.

## 4. Conclusion

In this paper, we have introduced the concept of bipolar neutrosophic soft expert set which is more effective and useful and studied some of its properties. Also the basic operations on neutrosophic soft expert set namely complement, union and intersection have been defined.